# A lower bound for Heilbronn's triangle-problem

G. Ellmann, 2025

methus02@hotmail.com

## Abstract

A conjecture by H. A. Heilbronn (before 1950) stated that if $n$ points are placed on the convex domain of unit area at a plane (no three points on a straight line), some of the triangles defined by these points always have an area not larger than $O(n^{-2})$, irrespective to the placing of the points.

Komlós, Pintz and Szemerédi [***K***] disproved this in 1981 by showing that $n$ points can be placed on the convex domain of unit area so that any of the triangle's area is not less than $O(n^{-2} \log n)$.

In this paper we are going to use a heuristic method to prove the following:

***Theorem 1***: *There is a constant $n_1 \in \mathbb{R}_+$ and a construction containing a Set of $n_0$ points within the unit circle of the plane where all the triangles defined by the points of this Set have an area $A_u \geq c\, n_0^{-3/2} log^{-7/2} n_0$ if $n_0 > n_1$.* Here $c \leq 1.655$…

The unit circle has an area of $\pi$. To satisfy the prescription of the original conjecture, i.e., the convex domain of area 1, one shall change the value of the constant factor accordingly.

## Introduction

### Conventions

All quantities in this paper are real.

Small Latin letters are used for positive integers, except for $e$, the base of the natural logarithm. Another exception is $c$, which denotes a positive constant, not necessarily the same from case to case. The letters $p$, $q$ and $r$ are reserved to denote radii of the circles examined. The expressions of $p$-circle and $p$-polygon are often used. These mean a circle with radius $p$, and a regular $p$-sided polygon written into that circle, respectively.

Capital letters denote points, straight lines, sections or real numbers, except $N$, which is an integer. The letter $A$ is reserved for the triangle's area on the disk, of course, it is variable from case to case. Sometimes it is denoted with $A_d$, where $d$ refers to the disk, while $A_u$ refers to the unit circle. $H$ is the height over the base side of a triangle. $P$, $Q$ and $R$ denote some elements of the Set. $O$ stands for the Origo, the centre of the circles. O() is the Ordo in the sense of Bachmann - Landau, i.e. the magnitudes of the remnants that are added to the main term, getting the correct value of the expression.

The Greek letters $\lambda$ and $\mu$ are used for angles, the $\Delta$ means triangle, while $\pi$ retains its well-known numeric value. As usual, $\pi()$ is standing for the number of primes from 2 to the argument. The # symbol is used for cardinality, i.e. it denotes the "number of".

The disk holds all the concentric circles we have, but none of smaller or larger circles. The Set contains all the points we are working with.

A straight line is denoted by naming two points on it, e.g., *PQ* line. The word "line" is always there. Sometimes a single capital letter denotes the line. The distance between points *P* and *Q* is denoted by *PQ*, in that case the word "line" is omitted.

**The initial Set**

Let $N$ and $p$ be primes, with $p$ in the interval $(N - [N/\log N], N - 2)$. Take the $p$-sided regular polygon on the plane, written into a circle with radius $p$ around the Origo. The word "circle", from now on, means one of these concentric circles, unless it is noted to be different from them. The $p$ prime run through the above interval, generating several polygons. Let a vertex of each polygon be on the $+X$ axis.

Then we rotate the polygons around the Origo into anticlockwise, each of them with a different small angle $cs/N$, where $s \equiv t^2 \pmod N$ where $t$ is the radius of the circle containing the polygon and $s$ has the smallest positive number satisfying the congruence. Here $0 < c$ is a constant, having the same value for all polygons. Every vertex of all the rotated polygons are the initial elements of our Set.

## ***Theorem 2***: *All the polygons are rotated with different angles.*

**Proof of *Theorem 2*:** We may assume that two polygons are on different circles. These circles have radii $p < q$.

Clearly, if $u \equiv p^2$ and $v \equiv q^2 \pmod N$ then $u \neq v$ since $|v - u| \equiv |q^2 - p^2| = (q - p) \cdot (q + p) \not\equiv 0 \pmod N$ because $0 < q - p < N$ and $N < q + p < 2N$ due to $\min(p, q) > ½N$ and $\max(p, q) < N$. Thus, the angles of rotations are all different. Therefore, they are still different when we multiply them with any value for $c$, the same factor for all the polygons.

By this ***Theorem 2*** are proven.

The side effect of $c$ in the rotation is that $c$ can be chosen so, that there are no three points on the Set that are collateral, and no *PQ* line that is tangent to a circle. If there is a *P*, *Q* pair that is collateral with some *R* of the Set or tangent to a circle after the rotation, then there is finite number of them. But $c$ is an element of a set with infinite cardinality, so there is always selectable infinite $0 < c < 1$, causes that the above situation cannot happen. Choose such a $c$ value!

**The cardinality of the initial Set**

According to Pierre Dusart [***D***], simplified a bit, $\pi(x) = x/\log x + x/\log^2 x + O(x/\log^3 x)$, if $x$ is large.

Using the above, it is easy to verify that $\pi(N(1 - 1/\log N)) = N/\log N + O(N/\log^3 N)$, if $N$ is large enough. Therefore, the number of polygons, i.e., the number of primes in the interval is $\pi(N - 2) - \pi(N - N/\log N) = N/\log^2 N + O(N/\log^3 N)$. Each polygon has at least $[N - N/\log N]$ vertices and at most $N - 2$ vertices, therefore the cardinality of the initial Set is at least and at most $n = [N^2/\log^2 N + O(N^2/\log^3 N)]$ where the Ordo bears different constants for the upper and lower estimates. Thus, the cardinality of the initial Set is $n \approx N^2/\log^2 N$.

**Choice of the names for the vertices**

If the three vertices of the Δ is on the same polygon then $A_d = r \sin(2\pi/r) \cdot r(1 - \cos(2\pi/r) = (2\pi + O(r^{-2})) \cdot (r(½(2\pi/r)^2 + O(r^{-4})) \geq 4\pi^3 N^{-1}\log N + O(N^{-3}\log^3 N)$. We can do this $A_d$ to be a bit smaller while we remain on the disk, but the three vertices are not on a single standard polygon anymore.

Let us assume that they are placed on two or three polygons. The radii of the circles on what *P*, *Q* and *R* vertices reside are called *p*, *q* and *r*. Let *PQ* be the "base" of the *PQR* Δ.

Observe that the *PQ* line cuts the *r* circle. It is because if the *PQ* line does not cut the *r* circle, then we rename the $r < p < q$ circles to $p < q < r$.

**The hit points and the exclusion zones**

The intersection points of the *PQ* lines with the *r* circle are referred to as **hit points**. In fact, we assume that the hit points are uniformly distributed on the *r* circle, also between neighbouring *R* vertices. This assumption makes the article heuristic.

The 2*Z* long interval around an *R* vertex is called an **exclusion zone**, or just **zone**, because we might exclude from the Set the *R* vertex where the hit-point is falling into that zone for some *P* and *Q* points. The exclusion zone is served to ignore the *PQR* Δs that have too small an area. We will determine max(*Z*) later.

It is important to remark that we cannot say that a certain point will be kept all the time as an element of the final Set because we don't know which points will be dropped. Some of the points that are retained in a certain phase might be eliminated in a later phase of the dropping procedure. We can, however, say that we drop many of the initial points but not all. The final cardinality of the Set is relevant, not the elements that are retained.

**On the assumed uniform distribution of the hit points**

Although the Set is deterministic, nevertheless, when $N \to \infty$ the hit points around an *R* vertex can be considered random with assumed uniform distribution.

The *r* can be any of the radii (of course, prime), thus the *PQ* lines do not "know" the location of the *r* circles because large primes distance from each other are not predictable. It seems to be no reason to have more or fewer hit points on one arc of the *r* circle than another arc of the same size on the same circle. Therefore, we may assume that the hit points on the *r* circle are uniformly distributed. Although we cannot prove that the hit points (of any possible *P* and *Q*) are uniformly distributed in the $2\pi$ long arc between neighbouring *R* vertices, it may still be assumed.

The exclusion zone's concept is based upon the assumed uniform distribution of the *PQ* lines' hit points on the *r* circles. Since the zone's length is 2*Z*, a single *PQ* line's hit point has $Z/\pi$ probability to fall into the zone around the nearest *r*-vertex.

## Prove *Theorem 1*

To prove ***Theorem 1***, we need two more theorems and two ***Lemmas***.

## Hit points outside the zone

Of course, if *R* is far from the *PQ* line, then *H* (the height of the *PQR* Δ) is large, therefore *A* is great, hence *R* is uninteresting. Let *W* be a hit point on the *r* circle. We would examine the case when *R* is close to *PQ*, i.e., $RW < N^{-1/2}$.

***Theorem 3***: *Let Z be a parameter and let W be a hit point of the PQ line and the r circle, where the P and Q points (both belong to the Set) are on the p and q circles, respectively (p < q are primes). If the W point has the distance from the nearest R vertex on the r circle is greater than some Z* ($\le N^{-1/2}$) *then the area A of the PQR Δ is A > Z.*

**Proof of *Theorem3***

Let us name the circles on which the Δ's vertices reside, $p < q \le r$. Look at the $R$ points to be retained! Let $R''$ be on a straight line of $OR$.

$W$ is on the $r$ circle, while $R''$ is on the $T$ line, perpendicular to the $WO$ line. If $N$ is great and $RW < N^{-1/2}$ then the distance of $R''$ from the $r$ circle is very short ($RR'' \approx (r^2 + (N^{-1/2})^2)^{1/2} - r <$ ½$N^{-2}\log N + O(N^{-5}\log^3 N)$), so we make a negligible small error by using $R''$ instead of $R$, when $N$ is large. If $RW \ge N^{-1/2}$ then we retain $R$ in the Set.

There are three cases: $p < q < r$, $p < q = r$ and $r = p < q$, because the $r < p < q$ case is equal to the $p < q < r$ case by changed names.

For the retained $R$ point $Z = R''S < R''W \approx RW$. See ***Diagram 1*** below, showing the bulky $p < q < r$ case.

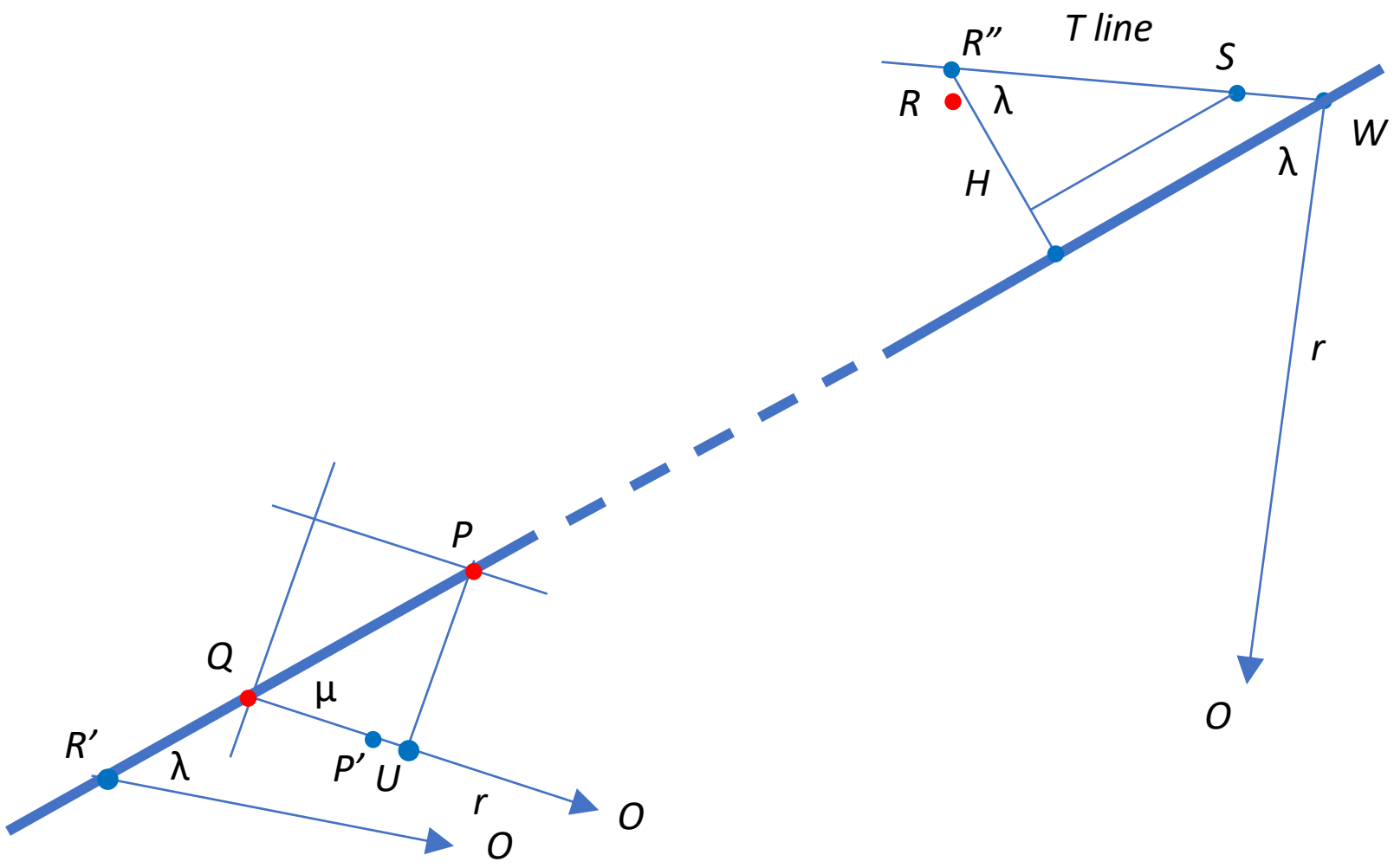

Let λ denote the smallest positive angle between the $PQ$ line and the $OW$ line. The $r$ circle's and the $PQ$ line's shared point that is other than $W$ is denoted by $R'$. Thus, $OR' = OW$ and $H \approx R''W \cdot \cos\lambda \approx RW \cdot \cos\lambda > Z\cos\lambda$.

$U$ is the projection of $P$ onto the $OQ$ line. The $p$ circle cuts the $OQ$ line at the $P'$ point. The $q$ circle is greater than the $p$ circle, so the tangent at $Q$ is further away from $P$ than $\min(p - q) = 2$ (the distance is $P'Q$). A side effect of $c$ in the rotation of the polygons (*Theorem 2*) ensures that $PQ$ cannot be tangent to any circle. From this follows that λ cannot be $\pi/2$, so $\cos\lambda \ne 0$. Thus, $P'Q = PQ\cos\mu \ge 2$. Because $r \ge q$, then $\mu \ge \lambda$. Since we are in the $<0, \pi/2>$ interval, $\cos\mu \le \cos\lambda$, therefore $PQ\cos\lambda \ge 2$. We have seen that $H > Z\cos\lambda$. Finally, $A =$ ½$PQ \cdot H >$ ½ $\cdot 2/\cos\lambda \cdot Z\cos\lambda$. Therefore ***A > Z***. By this, ***Theorem 3*** is proven for the $p < q < r$ case.

The $p < q = r$ case is essentially not deviate from the $p < q < r$ case, so, the proof is the same.

For the $r = p < q$ case: Here is $\cos\mu > \cos\lambda$, hence this case needs special treatment. For a moment, we name the circles to $p = q < r$. $RW > Z$, therefore $H \approx RW \cdot \cos\lambda > Z\cos\lambda$, as before. Now, $\cos\mu < \cos\lambda$, because $q < r$, so, $PQ\cos\lambda \ge PQ\cos\mu = QU$. But $QU$ runs from $cp\sin(\pi/p) = c\pi + O(N^{-2}\log^2 N)$ to about $cN/\log N$ and falls back to $cp\sin(2\pi/p) = 2c\pi + O(N^{-2}\log^2 N)$, because $P$ and $Q$ are on the same circle. Thus, in all circumstances $QU > 2$, hence $A >$ ½ $\cdot 2/\cos\lambda \cdot Z\cos\lambda$. To finish the proof of this case we reset the original names.

By this, ***Theorem 3*** is proven for all the cases.

## How large can the *Z* parameter (threshold) be?

### The principles

When the hit point has a distance from *R* that is greater than *Z*, then *R* is retained in the Set, otherwise the *PQR* Δ might be eliminated by eventually removing *R* from the Set. *Z* shall be chosen as large as possible, but small enough that it not causes to drop the whole Set.

We apply the following strategy: Fix the *Q* point, temporarily. Now we look at the *PQ* lines for all possible *P*. We collect the *PQ* bases into two classes (the first of these has several subclasses according to the increasing length of *PQ*). Each subclass will limit the #*R*, i.e., the number of *r*-vertices that can be near to the given *PQ* line. The heuristic assumption was that the probability of the hit point falling into the exclusion zone around the nearest *R* is less than $Z/\pi$, so we can estimate the number of points that might be eliminated.

***Theorem 4***: *If the Δ's three vertices are elements of the Set and they fall on two or three polygons, then max(Z) =* $c\,N^{-1}/\sqrt{\log N}$.

### Proof of *Theorem 4*

Let *Q* be any vertex of the Set. Thus $\#Q \le N^2/\log^2 N$, i.e., #*Q* is (over)estimated to be at most as big as the Set's initial cardinality. Let us estimate #*P* and by that the number of *PQ* lines. We use this moment also to estimate #*R* for the *R* points near that *PQ* line. After those estimates, we let *Q* run through all the points of the Set.

For a fixed *Q*, the #*P* can be (over)estimated, based on the half area of a *U* circle (a special circle, not concentric with the other circles), centred in *Q* and having radius of max(*PQ*). The reason for using half of the circle is $p < q$.

We create many subclasses for the length of the bases, namely $\frac{1}{4}e^{1+j} \le P_jQ < \frac{1}{4}e^{2+3j/2}$
$(1 \le j \le [\frac{2}{3} \cdot \log N + 19/18])$

The additive 19/18 ensures that the end of the last interval is always greater than 2*N*, the diameter of the largest circle. So, we have at most $\frac{2}{3} \cdot \log N + 19/18$ subclasses. The intervals overlap; each interval starts before the previous one ends. When $P_jQ$ falls into both intervals, we count it to the smallest index.

Let us estimate $\#P_j$.

***Lemma 1***: $\#P_j < \frac{1}{2}\pi\, e^{2+j}$.

### Proof of *Lemma 1*

Denote $G_j = \min(P_jQ)$ within a subclass. Consider a grid made of quadrats, each having an edge $G_j$ long. The #nodes of the grid positioned in the half of circle *V* around centre *Q* with radius max(*PQ*) is an (over)estimate for #*P* belongs to a fixed *Q*. If there were more points to be counted, those would be placed on the board of some squares or inside them, hence, two points would be closer to each other than $G_j$. This would be a contradiction with $G_j = \min(P_jQ)$. Former subclasses would have dealt with this. So, the area of the half circle *V* is divided by the $G_j^2$ area of the smallest square in the grid, and this could give us an over-estimate for the number of nodes of the grid. (The node at the upper right would stand for the square.)

Note, that there is well-known waste at the edge of the half circle. It is at most $O((\#P_j)^{1/2})$ according to C. F. Gauss (Circle-problem) [***G***]. This waste is marginal; thus, we can disregard it.

$\#P_j < ½\pi\ (\max(P_jQ))^2/(\min(P_jQ))^2 = ½\pi\ e^{4+3j}/e^{2+2j} = ½\pi\ e^{2+j}$.

By this ***Lemma 1*** is proven.

Now let us estimate $\#R_j$.

***Lemma 2:*** $\#R_j \le 8\sqrt{2} \cdot N/\sqrt{\log N} \cdot e^{-1-j}$.

**Proof of *Lemma 2***

We look at the $R$ points near to the $P_jQ$ line. We know that their number is at most $L/G_j$ in this subclass. This is so, because two $R$ points that are nearer to each other than $G_j$, would belong to previous subclasses as $P$s.

The length of $L$ is an overestimate of the useful part of the $P_jQ$ line that consists of one or two sections. $R_j$ cannot be outside the largest circle, nor can it be closer to the origo than the smallest circle's radius. Hence, the length of $L$, or the joint length of the two sections will not exceed the length of the longest continuous section, i.e., twice the length of the tangent from a point on the outermost circle drawn to the innermost circle (see the Appendix). So, $L$ is $2\sqrt{2} \cdot N/\sqrt{\log N} + O(N/\log^{-3/2}N)$. Ultimately, $\#R_j \le 8\sqrt{2} \cdot N/\sqrt{\log N} \cdot e^{-1-j}$.

By this, ***Lemma 2*** is proven.

The overall number of Δs with base of $P_jQ$ on the $P_jQ$ line which have very small area for this subclass is $\#P_j \cdot \#Q \cdot \#R_j < (½\pi\ e^{2+j}) \cdot (N^2/\log^2 N) \cdot (8\sqrt{2}\ N/\sqrt{\log N} \cdot e^{-1-j}) = 4\pi e\sqrt{2}\ N^3/\log^{5/2}N$.

How many Δs might be rejected due to their area too small? There is at most $\#P_j \cdot \#Q \cdot \#R_j \le \pi e \cdot 4\sqrt{2} \cdot N^3/\log^{5/2}N$ points as candidate for deletion in this subclass of $P_j$. Since there is no $j$ on the right-hand-side, we might assume that the same $Z$ serves for all the subclasses. Overall, we might drop ¾ $N^2/\log^2 N$ points from the Set by treating each of the [⅔ · log $N$ + 19/18] subclasses, so there are at most $2\pi e \cdot 4\sqrt{2}/3 \cdot N^3/\log^{3/2}N + O(N^3/\log^{5/2}N)$ points that will be deleted with $Z/\pi$ probability. Therefore, we drop at most $8\pi e\sqrt{2}/3 \cdot N^3/\log^{3/2}N \cdot Z/\pi = ¾\ N^2/\log^2 N$ points. (The ¾ here can be changed to any other number < 1 because we do not want to eliminate all the points. We choose ¾ for convenience.) Thus, $\mathbf{\max(Z) = c_2 N^{-1}/\sqrt{\log N}}$ where $c_2 = 9\sqrt{2}/(4e) = 1.14\ldots$ .

Hence, for the vertices resting on two or three polygons, $Z = c_2 N^{-1}/\sqrt{\log N}$. By this, ***Theorem 4*** is proven.

## Summary

After dropping at most ¾ $N^2/\log^2 N$ points of the original cardinality, the Set shall retain at least a quarter of its original points.

In the following the subscript of $Z$ refers to the number of circles that give home to the $P$, $Q$ and $R$ vertices of the Δ.

Trivially, $\max(Z_1) \approx 4\pi^3 N^{-1}\log N$ far too large, due to $A \ge 4\pi^3 N^{-1}\log N + O(N^{-3}\log^3 N)$. The consequence of ***Theorem 4*** is $\max(Z_2) = \max(Z_3) = c_2 N^{-1}/\sqrt{\log N}$. From these we conclude that $Z = c_2 N^{-1}/\sqrt{\log N}$, if $N$ is great enough. So, after the deletions, constructed from the remaining points the Δs area $A$ is certainly greater than $\max(Z) = c_2 N^{-1}/\sqrt{\log N}$.

## Mapping onto the unit circle as domain

The Set's starting cardinality was $n \approx N^2/\log^2 N$. After the drop this became $n \geq ¼\, N^2/\log^2 N$. One can verify that the inverse of this final cardinality is $N \leq ¼ n^{1/2} \log n + O(n^{1/2} \log^{-1} n)$. Please see that $\log N = ½ \log n + O(\log\log n)$.

The disk-ring domain has been used throughout this paper. Now we project our result onto the unit circle, i.e., one that has radius 1. We have seen that for the Δ having the smallest area $A$ constructed from the remaining points of the Set has $A_d \geq \max(Z) = c_2\, N^{-1}/\sqrt{\log N} \approx 2c_2\, n^{-1/2}\log^{-3/2} n$. For the mapping we must divide the disk's $A_d$ by $N^2 \leq n \log^2 n + O(n)$. Therefore, on the unit circle, $A_u \geq c\, n^{-3/2}\log^{-7/2} n + O(n^{-3/2} \log^{-11/2} n)$, if $n$ is large enough. Here $c \leq c_2\sqrt{2} = 9/(2e) = 1.655\ldots$

By this, the ***Theorem 1*** is proven.

## Acknowledgement

The author is in debt to Prof. E. Szemerédi for his suggestions and advises.

## Appendix: on *L* and *DS*

Let $a, b, X \in \mathbb{R}_+$, $a < b$, and let $X \leq a$. See ***Diagram 2*** below.

Trivially, $a^2 - X^2 < b^2 - X^2$. It is also clear that $a^2 - X^2 < \sqrt{(b^2 - X^2)(a^2 - X^2)}$. From this follows $a^2 + b^2 - 2X^2 - 2\sqrt{(a^2 - X^2)(b^2 - X^2)} \leq b^2 - a^2$. Thus, taking the meaningful positive square root, $\sqrt{b^2 - X^2} - \sqrt{a^2 - X^2} \leq \sqrt{b^2 - a^2}$. This means that $FM \leq JK = L$ when $0 < X \leq a$. If $a < X \leq b$, then $JK = L$ is obviously greater than $BE$.

For those who needs the solution of $b - a \leq \sqrt{b^2 - X^2} - \sqrt{a^2 - X^2}$ without differentiation, here is the summary of that: the $-2abX^2 \geq -X^2(a^2 + b^2)$ leads to $ab - X^2 \geq \sqrt{(b^2 - X^2)(a^2 - X^2)}$ which is equivalent to the searched inequality, taking again the meaningful positive square root.

***Diagram 2***

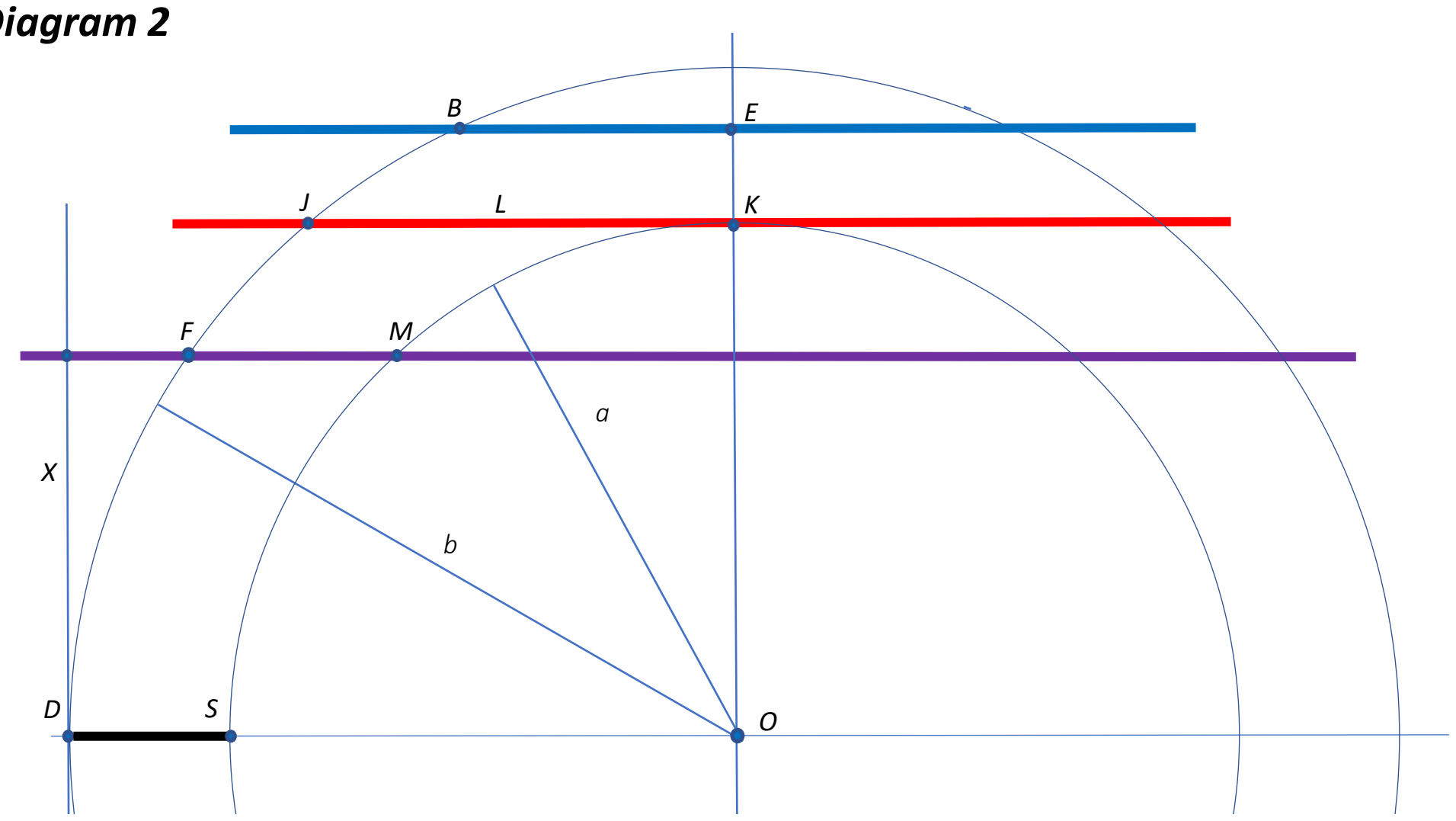